\documentclass[pdflatex,sn-mathphys-num]{sn-jnl}
\usepackage[utf8]{inputenc}
\usepackage{amsmath,amssymb,amsfonts,mathtools}
\allowdisplaybreaks
\theoremstyle{thmstyleone}
\newtheorem{theorem}{Theorem}[section]
\newtheorem{lemma}[theorem]{Lemma}

\theoremstyle{thmstylethree}
\newtheorem{problem}{Problem}
\newtheorem{remark}{Remark}
\newcommand{\Sym}{\operatorname{Sym}}
\newcommand{\Aut}{\operatorname{Aut}}
\newcommand{\AGL}{\operatorname{AGL}}
\newcommand{\Stab}{\operatorname{Stab}}
\newcommand{\diag}{\operatorname{diag}}
\newcommand{\Cay}{\operatorname{Cay}}
\newcommand{\GL}{\operatorname{GL}}
\newcommand{\F}{\mathbb F}
\newcommand{\Q}{\mathbb Q}
\newcommand{\Z}{\mathbb Z}
\newcommand{\calB}{\mathcal B}
\newcommand{\calE}{\mathcal E}
\newcommand{\calC}{\mathcal C}

\newcommand{\id}{\mathrm{id}}

\newcommand{\supp}{\operatorname{supp}}

\renewcommand{\theenumi}{\roman{enumi}}

\hypersetup{pdftitle={A Sylow reduction for isomorphisms of Cayley graphs on abelian groups of cube-free order},pdfauthor={Tingzeng Wu, Liang Qin, Yong Yu and Xing Gao}}

\begin{document}

\title{A Sylow reduction for isomorphisms of Cayley graphs on abelian groups of cube-free order}

\author*[1,2]{\fnm{Tingzeng} \sur{Wu}}\email{mathtzwu@163.com}
\author[1]{\fnm{Liang} \sur{Qin}}\email{ql22120806@126.com}
\author[3]{\fnm{Yong} \sur{Yu}}\email{yuyong\_lzu@yeah.net}
\author[3,4]{\fnm{Xing} \sur{Gao}}\email{gaoxing@lzu.edu.cn}

\affil[1]{\orgdiv{School of Mathematics and Statistics},
  \orgname{Qinghai Minzu University},
  \orgaddress{\city{Xining}, \state{Qinghai} \postcode{810007}, \country{P.R.~China}}}
\affil[2]{\orgname{Qinghai Institute of Applied Mathematics},
  \orgaddress{\city{Xining}, \state{Qinghai} \postcode{810007}, \country{P.R.~China}}}
\affil[3]{\orgdiv{School of Mathematics and Statistics, Key Laboratory of Applied Mathematics and Complex Systems},
  \orgname{Lanzhou University},
  \orgaddress{\city{Lanzhou}, \postcode{730000}, \country{P.R.~China}}}
\affil[4]{\orgname{Gansu Provincial Research Center for Basic Disciplines of Mathematics and Statistics},
  \orgaddress{\city{Lanzhou}, \postcode{730070}, \country{P.R.~China}}}

\abstract{The isomorphism problem for Cayley graphs is one of the fundamental problems in algebraic combinatorics. Dobson studied this isomorphism problem and asked whether, assuming that the group generated by the left regular representation and its conjugate under an isomorphism admits a complete sequence of normal block systems, the isomorphism problem for Cayley digraphs on a finite abelian group of cube-free order can be reduced to the corresponding problem on its Sylow subgroups. We answer this question affirmatively. We first prove that color-preserving isomorphisms preserving the fiber partition on direct products $A\times H$ of abelian groups of coprime orders, where $|H|$ is a prime or the square of a prime, have a product decomposition property. Applying this property successively to the Sylow subgroups of a cube-free-order abelian group $G$, we show that if two colored binary Cayley structures on $G$ are connected by a color-preserving isomorphism $f$ for which $\langle G_L,f^{-1}G_Lf\rangle$ admits a complete sequence of block systems, then they are also connected by a color-preserving Sylow coordinate isomorphism. This yields the corresponding Sylow reduction for solving sets.}


\maketitle

\section{Introduction}\label{sec: introduction}

Let $G$ be a finite group and let $S\subseteq G\setminus\{1\}$. The Cayley digraph $\Cay(G,S)$ has vertex set $G$, connection set $S$, and arc set $\{(g,gs):g\in G,\ s\in S\}$. When $S=S^{-1}$, we regard it as an undirected graph. A bijection $f:G\to G$ is an isomorphism from $\Cay(G,S)$ to $\Cay(G,T)$ if and only if
\[
 g^{-1}h\in S\quad\Longleftrightarrow\quad f(g)^{-1}f(h)\in T
 \qquad(g,h\in G).
\]

A group $G$ is a \emph{DCI-group} (directed Cayley isomorphism group) if, whenever two Cayley digraphs on $G$ are isomorphic, some automorphism of $G$ maps one connection set onto the other. Let $G=G_1\times\cdots\times G_r$ be a finite abelian group written as the direct product of its Sylow subgroups. Let $\Sym(G)$ denote the symmetric group on $G$, and let $D_G=\prod_i\Sym(G_i)\leq\Sym(G)$, where $(h_i)$ acts by $(x_1,\dots,x_r)\mapsto(h_1(x_1),\dots,h_r(x_r))$. We call an isomorphism in $D_G$ a \emph{Sylow coordinate isomorphism}. We also let $G_p$ denote the Sylow $p$-subgroup of $G$. We say that the isomorphism problem reduces to the Sylow subgroups if any two isomorphic Cayley structures on $G$ admit a Sylow coordinate isomorphism between them. Let $C_m$ denote the cyclic group of order $m$ and $C_m^k$ the direct product of $k$ copies of $C_m$. The left regular representation of $G$ is $G_L=\{L_g:g\in G\}$, where $L_g(x)=gx$. For any group $A$, we let $A_L$ denote its left regular representation.

The isomorphism problem for Cayley graphs is one of the classical problems in algebraic combinatorics~\cite{li2002}. Muzychuk~\cite{muzychuk1995,muzychuk2004} proved that cyclic groups of square-free order are DCI-groups and solved the isomorphism problem for colored Cayley graphs on cyclic groups of arbitrary order. Here an integer is square-free if each of its prime divisors occurs with multiplicity one. For a bijection $f\colon G\to G$, Dobson~\cite{dobson2014} proved that if the group generated by two regular abelian permutation groups $G_L$ and $f^{-1}G_Lf$ admits a complete sequence of block systems, then one of the two groups can be conjugated within this group so that the resulting pair generates a solvable group.
Feng and Kov\'acs~\cite{feng2018} proved that the elementary abelian group $C_p^5$ is a DCI-group for every odd prime $p$. For an abelian group $G$ of order $n$, Dobson~\cite{dobson2018} subsequently established a reduction of solving sets to the Sylow subgroups under the assumptions that
\begin{equation}\label{gcd}
\gcd\bigl(\operatorname{rad}(n),\varphi(\operatorname{rad}(n))\bigr)=1
\end{equation}
and that, for every $f\in\Sym(G)$, the group generated by $G_L$ and $f^{-1}G_Lf$ admits a complete sequence of normal block systems. Here $\operatorname{rad}(n)$ is the product of the distinct prime divisors of $n$, and $\varphi$ is Euler's totient function. Kov\'acs et al.~\cite{kovacs2023} proved that $C_p\times C_n$ ($p$ prime, $n$ square-free) and $C_p^2\times C_q^2$ ($p,q$ distinct primes) are also DCI-groups. These results show that the Sylow reduction method is an effective approach to the isomorphism problem for Cayley graphs.

An integer $n$ is cube-free if $p^3\nmid n$ for every prime $p$, and an abelian group is of cube-free order if its order is cube-free. Such orders need not satisfy~\eqref{gcd}, which led Dobson to pose the following problem.

\begin{problem}[{Dobson~\cite{dobson2018}}]\label{prob: dobson}
Suppose that the group generated by $G_L$ and $f^{-1}G_Lf$ admits a complete sequence of normal block systems. Does the isomorphism problem for Cayley digraphs on an abelian group of cube-free order reduce to the corresponding problem for its Sylow subgroups?
\end{problem}

Since $G_L$ is a regular abelian subgroup, the following text will show that every invariant block system is normal, hence the normal and invariant formulations of the block-chain hypothesis agree.

We answer Problem~\ref{prob: dobson} affirmatively. To state our results, we first introduce some terminology. Let $\Omega$ be a set of size $n$, and let $K\leq\Sym(\Omega)$ be transitive. A partition $\calB$ of $\Omega$ is a collection of pairwise disjoint nonempty subsets, called \emph{blocks}, whose union is $\Omega$. We say that $\calB$ is \emph{$K$-invariant}, or a \emph{block system} for $K$, if $g(B)\in\calB$ for all $g\in K$ and $B\in\calB$. For a binary relation $Q\subseteq\Omega^2$ and $g\in K$, let $g(Q)=\{(g(x),g(y)):(x,y)\in Q\}$. The relation $Q$ is $K$-invariant if $g(Q)=Q$ for all $g\in K$.

Let $t=\Omega(n)=\sum_{p\mid n}v_p(n)$, the number of prime factors of $n$, counted with multiplicity. A chain of $K$-invariant partitions $\calB_0\prec\calB_1\prec\cdots\prec\calB_t$ is a \emph{complete sequence of block systems} for $K$ if it starts with the discrete partition $\{\{x\}:x\in\Omega\}$ and ends with the trivial partition $\{\Omega\}$, and the ratio of successive block sizes is prime at every step. The sequence is \emph{normal} if each partition is the orbit partition of a normal subgroup of $K$.

A \emph{colored binary Cayley structure} on $G$ is a pair $X=(G,(R_c)_{c\in\calC})$, where $\calC$ is a finite set of colors, $R_c\subseteq G^2$ for each $c\in\calC$, and $G_L\leq\Aut(X)$. Here $\Aut(X)$ consists of the permutations of $G$ preserving each relation $R_c$. For another such structure $Y=(G,(R'_c)_{c\in\calC})$ with the same color set, a bijection $f:G\to G$ is a \emph{color-preserving isomorphism} from $X$ to $Y$ if, for every $c\in\calC$ and all $g,h\in G$, $(g,h)\in R_c$ exactly when $(f(g),f(h))\in R'_c$.

Let $\pi_H:A\times H\to H$ be the projection $(a,h)\mapsto h$. For $h\in H$, we call $\pi_H^{-1}(h)=A\times\{h\}$ the \emph{$A$-fiber over $h$}, and we call $\calE_A=\{A\times\{h\}:h\in H\}$ the \emph{$A$-fiber partition}. A bijection $f:A\times H\to A\times H$ preserves the fiber partition if it maps each $A$-fiber onto an $A$-fiber.

For direct products $A\times H$ of abelian groups of coprime orders, where $|H|$ is a prime or the square of a prime, we first prove a product decomposition property for color-preserving isomorphisms that preserve the $A$-fiber partition.

\begin{theorem}\label{1.1}
Let $A$ be a finite abelian group, let $q$ be a prime not dividing $|A|$, and let $H$ be an abelian group of order $q$ or $q^2$. Let $X$ and $Y$ be colored binary Cayley structures on $A\times H$, and suppose that
\[
 f(a,h)=(\beta_h(a),\tau(h))\colon X\longrightarrow Y
\]
is a color-preserving isomorphism that preserves the fiber partition. Put $R=(A\times H)_L$ and $B=\langle R,fRf^{-1}\rangle\leq\Aut(Y)$. Then there exists $b\in B$ such that, on writing $bf(a,h)=(\widetilde\beta_h(a),\widetilde\tau(h))$, the map $(a,h)\mapsto(\widetilde\beta_{h_0}(a),\widetilde\tau(h))$ from $X$ to $Y$ is a color-preserving isomorphism for every $h_0\in H$.
\end{theorem}

Applying the above product decomposition property to the Sylow subgroups of a finite abelian group of cube-free order yields the following result.

\begin{theorem}\label{thm: main}
Let $G$ be a finite abelian group of cube-free order, and let $X$ and $Y$ be colored binary Cayley structures on $G$. Suppose that there is a color-preserving isomorphism $f:X\to Y$ such that $\langle G_L,f^{-1}G_Lf\rangle$ admits a complete sequence of block systems. Then there is a color-preserving isomorphism $d=\prod_{p\mid |G|}d_p:X\to Y$ with $d_p\in\Sym(G_p)$ for each prime $p$ dividing $|G|$.
\end{theorem}

Let $X$ be a colored binary Cayley structure on $G$. A set $S\subseteq\Sym(G)$ is a \emph{solving set} for $X$ if, for every colored binary Cayley structure $Y$ on $G$ isomorphic to $X$, some $s\in S$ is a color-preserving isomorphism from $X$ to $Y$. If $S$ is a solving set for every finite colored binary Cayley structure on $G$, we call it a \emph{universal color-preserving solving set} for $G$. The following result gives the corresponding Sylow reduction for solving sets.

\begin{theorem}\label{thm: solvingset}
Let $G$ be a finite abelian group of cube-free order and let $X$ be a Cayley digraph on $G$. Suppose that, for every Cayley digraph $Y$ on $G$ isomorphic to $X$, there is an isomorphism $f:X\to Y$ for which $\langle G_L,f^{-1}G_Lf\rangle$ admits a complete sequence of block systems. For each $p\mid |G|$, let $S_p$ be any universal color-preserving solving set for $G_p$. Then $\prod_{p\mid |G|}S_p$ is a solving set for $X$.
\end{theorem}

\begin{remark}\label{rem: problem1}
Theorems~\ref{thm: main} and~\ref{thm: solvingset} jointly give an affirmative answer to Problem~\ref{prob: dobson}.
\end{remark}

The remainder of this paper is organized as follows. Section~\ref{sec: prelim} collects the preliminary results. We prove Theorem~\ref{1.1} in Section~\ref{sec: fourier}; Section~\ref{sec: proof} proves Theorems~\ref{thm: main} and~\ref{thm: solvingset}. The final section concludes the paper.

\section{Preliminaries}\label{sec: prelim}

This section collects some preliminaries needed in the sequel. Throughout, all groups and sets are finite unless otherwise stated. We compose permutations from right to left, so $(ga)(x)=g(a(x))$, and let $U^a=a^{-1}Ua$ denote the image of a subset $U$ under conjugation by $a$. For a subgroup $K\leq G$, let $N_G(K)=\{g\in G:gKg^{-1}=K\}$ denote its normalizer and let $C_G(K)=\{g\in G:gk=kg\text{ for all }k\in K\}$ denote its centralizer in $G$. Let $\pi(n)$ denote the set of prime divisors of a positive integer $n$.

For a set $\pi$ of primes, a finite group $G$ is a \emph{$\pi$-group} if $\pi(|G|)\subseteq\pi$. Let $R$ be a finite nilpotent group, let $R_p$ denote its Sylow $p$-subgroup, and let $R_{p'}=\prod_{q\ne p}R_q$ denote its Hall $p'$-subgroup, also called the $p$-complement.

Let $K\leq\Sym(\Omega)$. The point stabilizer of $x\in\Omega$ is $\Stab_K(x)=\{g\in K:g(x)=x\}$, and the setwise stabilizer of $D\subseteq\Omega$ is $\Stab_K(D)=\{g\in K:g(D)=D\}$. A permutation group is \emph{semiregular} if all its point stabilizers are trivial, and is \emph{regular} if it is both transitive and semiregular. For a $K$-invariant subset $B\subseteq\Omega$, let $K^B$ denote the permutation group induced by $K$ on $B$. For a $K$-invariant partition $\calB$, let $K^{\calB}$ denote the group induced on its blocks. The kernel of this action is $\operatorname{fix}_K(\calB)=\{g\in K:g(B)=B\text{ for every }B\in\calB\}$. A subgroup $V\leq U$ is \emph{characteristic} in $U$, written $V\operatorname{char}U$, if every automorphism of $U$ maps $V$ onto itself.

Let $C_p^2=C_p\times C_p$. Let $\F_q$ denote the finite field of order $q$, $\GL(n,q)$ the general linear group of degree $n$ over $\F_q$, and $\AGL(n,q)=\GL(n,q)\ltimes\F_q^n$ the corresponding affine general linear group. The exponent $\exp(A)$ of a finite group $A$ is the least common multiple of the orders of its elements. For a prime $q$ and a positive integer $m$, let $v_q(m)$ be the exponent of $q$ in the prime factorization of $m$, and let $\zeta_m$ denote a primitive $m$-th root of unity.

The support of a matrix $B=(b_{x,y})_{x,y\in\Lambda}$ indexed by a finite set $\Lambda$ is $\supp(B)=\{(x,y)\in\Lambda^2:b_{x,y}\ne0\}$. For a family $(B_c)_{c\in\calC}$ of such matrices, put $\supp((B_c)_{c\in\calC})=\bigcup_{c\in\calC}\supp(B_c)$. The digraph on $\Lambda$ whose arc set is this support is called the \emph{associated digraph} of the family. For a single matrix $B$, the digraph on $\Lambda$ with arc set $\supp(B)$ is called the associated digraph of $B$. A \emph{monomial matrix} has exactly one nonzero entry in each row and each column.

\begin{lemma}[Wielandt~\cite{wielandt1964}; orbit--stabilizer theorem]\label{std: orbitstab}
Let $K\leq\Sym(\Omega)$ and let $x\in\Omega$. The length of the $K$-orbit of $x$ is $|x^K|=[K:\Stab_K(x)]$.
\end{lemma}

\begin{lemma}[Hall~\cite{hall1928}; Hall's theorem]\label{std: hall}
Let $G$ be a finite solvable group, and let $\pi$ be a set of primes. Then the following hold.
\begin{enumerate}[\normalfont(iii)]
\item $G$ has a Hall $\pi$-subgroup.
\item Every $\pi$-subgroup of $G$ is contained in a Hall $\pi$-subgroup of $G$.
\item Any two Hall $\pi$-subgroups of $G$ are conjugate in $G$.
\end{enumerate}
\end{lemma}

\begin{lemma}[Sylow~\cite{sylow1872}; Sylow's theorem]\label{std: sylow}
Let $G$ be a finite group and let $q$ be a prime. Then the following hold.
\begin{enumerate}[\normalfont(iii)]
\item $G$ has a Sylow $q$-subgroup.
\item Every $q$-subgroup of $G$ is contained in a Sylow $q$-subgroup of $G$.
\item Any two Sylow $q$-subgroups of $G$ are conjugate in $G$.
\end{enumerate}
\end{lemma}

\begin{lemma}[Wielandt~\cite{wielandt1964}]\label{std: centralizer}
If $R\leq\Sym(\Omega)$ is regular and abelian, then $C_{\Sym(\Omega)}(R)=R$.
\end{lemma}

\begin{lemma}[Godsil and Royle~\cite{godsil2001}]\label{std: normalizer}
Let $G$ be a group and let $G_L$ be its left regular representation. Then
\[
N_{\Sym(G)}(G_L)=\{x\mapsto g\theta(x):g\in G,\ \theta\in\Aut(G)\}=\Aut(G)\ltimes G_L.
\]
In particular, if $G$ is abelian and written additively, then $\tau\in N_{\Sym(G)}(G_L)$ if and only if there are $g_*\in G$ and $\theta\in\Aut(G)$ such that $\tau(x)=g_*+\theta(x)$ for all $x\in G$. If $q$ is prime and $G$ is the additive group of $\F_q^n$, then $N_{\Sym(G)}(G_L)=\AGL(n,q)$.
\end{lemma}

\begin{lemma}[Robinson~\cite{robinson1996}]\label{std: abelian-sylow}
Let $R$ be a finite abelian group and let $S\leq R$. Then $S=\prod_{q\in\pi(|R|)}(S\cap R_q)$, where $R_q$ is the Sylow $q$-subgroup of $R$.
\end{lemma}

\begin{lemma}[Robinson~\cite{robinson1996}]\label{std: coprime-map}
Let $A$ be a finite abelian group and let $m$ be a positive integer coprime to $|A|$. Then the multiplication map $\mu_m:a\mapsto ma$ is a bijection on $A$. In particular, if $m\equiv1\pmod{\exp(A)}$, then $\mu_m=\id_A$. For any finite group $H$, if $m=|H|^t$ with $t\geq1$, then $h^m=1$ for every $h\in H$. In additive notation, this says that multiplication by $m$ is the zero map on $H$.
\end{lemma}

\begin{lemma}[Washington~\cite{washington1997}]\label{std: cyclotomic}
Let $m$ and $n$ be coprime positive integers. Then $\Q(\zeta_m,\zeta_n)=\Q(\zeta_{mn})$ and $[\Q(\zeta_m):\Q]=\varphi(m)$.
\end{lemma}

\begin{lemma}[Hardy and Wright~\cite{hardy2008}; Legendre's formula]\label{std: legendre}
Let $q$ be a prime and let $a$ be a positive integer. Then $v_q(a!)=\sum_{i\geq1}\lfloor a/q^i\rfloor$. In particular, $v_q((q^2)!)=q+1$.
\end{lemma}

In addition, we shall use the following Fourier diagonalization property. For a survey of the spectral theory of Cayley graphs, see Liu and Zhou~\cite{liu2022}.

\begin{lemma}[Terras~\cite{terras1999}]\label{std: fourier}
Let $A$ be a finite abelian group, let $\widehat A$ be its group of complex characters, and let $F=(\chi(a))_{a\in A,\ \chi\in\widehat A}$. Then the following hold.
\begin{enumerate}[\normalfont(iii)]
\item $\overline F^{\,T}F=|A|I$, so $F^{-1}=|A|^{-1}\overline F^{\,T}$.
\item $F$ simultaneously diagonalizes all $A$-circulant matrices, that is, all $|A|\times|A|$ matrices $M$ whose $(x,y)$-entry depends only on $y-x$. For each such $M$, the diagonal entry of $F^{-1}MF$ indexed by $\chi$ is the eigenvalue of $M$ corresponding to the eigenvector $\chi$.
\end{enumerate}
\end{lemma}

\begin{lemma}[{Dobson~\cite{dobson2014}}]\label{grp: normalblocks}
Let $K\leq\Sym(\Omega)$ be transitive and let $R\leq K$ be regular and abelian. Then every $K$-invariant partition is normal and is the orbit partition of a normal subgroup of $R$.
\end{lemma}

\begin{lemma}[{Dobson~\cite{dobson2014}}]\label{grp: solvable}
Let $n$ be a positive integer and put $\pi=\pi(n)$. Let $R$ and $T$ be regular abelian permutation groups on the same set $\Omega$ of size $n$. If $\langle R,T\rangle$ admits a complete sequence of block systems, then there is some $\delta\in\langle R,T\rangle$ for which $\langle R,T^\delta\rangle$ is a solvable $\pi$-group.
\end{lemma}

\begin{lemma}\label{grp: orbitbound}
Let $K\leq\Sym(\Omega)$ be transitive of degree $n=|\Omega|>1$, and let $\calB_0\prec\calB_1\prec\cdots\prec\calB_\ell$ be a complete sequence of block systems for $K$. Let $q=\min\pi(n)$ and let $n=q^as$ with $\gcd(q,s)=1$. If $H\leq K$ satisfies $\pi(|H|)\subseteq\{r\in\pi(n):r>q\}$, then $|x^H|\leq s$ for every $x\in\Omega$.
\end{lemma}

\begin{proof}
Fix $x\in\Omega$. For $0\leq i\leq\ell$, let $B_i(x)$ be the block of $\calB_i$ containing $x$, and put $H_i=\Stab_H(B_i(x))$. Each $\calB_i$ is $H$-invariant. Thus any element of $H$ stabilizing $B_{i-1}(x)$ also stabilizes the unique block $B_i(x)$ containing $B_{i-1}(x)$. It follows that $H_0\leq H_1\leq\cdots\leq H_\ell=H$. Since $H_0=\Stab_H(x)$, Lemma~\ref{std: orbitstab} gives
\begin{equation}\label{1}
|x^H|=[H:H_0]=\prod_{i=1}^{\ell}[H_i:H_{i-1}].
\end{equation}
Put $p_i=|B_i(x)|/|B_{i-1}(x)|$. Each $p_i$ is prime, and their product is $|B_\ell(x)|/|B_0(x)|=n$. The group $H_i$ induces a permutation group $L_i\leq\Sym(p_i)$ on the $p_i$ blocks of $\calB_{i-1}$ contained in $B_i(x)$. The stabilizer in $H_i$ of the block $B_{i-1}(x)$ is $H_{i-1}$, so $[H_i:H_{i-1}]\leq p_i$.

If $p_i=q$, then $|L_i|$ divides both $|H_i|$ and $q!$. Every prime divisor of $|H_i|$ is greater than $q$, since $H_i\leq H$ and $\pi(|H|)\subseteq\{r\in\pi(n):r>q\}$. Every prime divisor of $q!$ is at most $q$. Thus $L_i$ is trivial and $H_i$ fixes each block of $\calB_{i-1}$ contained in $B_i(x)$. In particular, $H_i=H_{i-1}$. Exactly $a$ of the factors in $\prod_i p_i=n=q^as$ are equal to $q$. Substituting the preceding bounds into~\eqref{1} therefore gives $|x^H|\leq\prod_{p_i\ne q}p_i=s$.
\end{proof}

\begin{lemma}\label{grp: alignment}
Let $R,T\leq\Sym(\Omega)$ be regular abelian groups on a set of size $n>1$. Suppose that $\calB_0\prec\calB_1\prec\cdots\prec\calB_\ell$ is a complete sequence of block systems for $\langle R,T\rangle$, and put $q=\min\pi(n)$. Then there exists $a_0\in\langle R,T\rangle$ such that $R_{q'}$ and $(T^{a_0})_{q'}$ have the same orbit partition $\calE$. Moreover, the following hold.
\begin{enumerate}[\normalfont(iii)]
\item $\calE$ is invariant under $\langle R,T^{a_0}\rangle$.
\item For each $B\in\calE$, the groups $R_{q'}^B$ and $(T^{a_0})_{q'}^B$ are regular on $B$.
\item For each $B\in\calE$, the group $\Stab_{\langle R,T^{a_0}\rangle}(B)^B$ admits a complete sequence of block systems on $B$.
\end{enumerate}
\end{lemma}

\begin{proof}
Let $n=q^as$ with $\gcd(q,s)=1$. By Lemma~\ref{grp: solvable}, there is some $b\in\langle R,T\rangle$ such that $K_1=\langle R,T^b\rangle$ is a solvable $\pi(n)$-group. Each $\calB_i$ is invariant under $K_1\leq\langle R,T\rangle$, so the given sequence is also a complete sequence of block systems for $K_1$.

Put $\pi_1=\pi(|K_1|)\setminus\{q\}$. The groups $R_{q'}$ and $(T^b)_{q'}$ are $\pi_1$-subgroups of $K_1$. By Lemma~\ref{std: hall}(ii), there are Hall $\pi_1$-subgroups $H_R$ and $H_T$ of $K_1$ containing $R_{q'}$ and $(T^b)_{q'}$, respectively. By part~(iii) of the same lemma, there is some $c\in K_1$ with $H_T^c=H_R$. Let $H=H_R$. Since $(T^{bc})_{q'}=(T^b)_{q'}^c\leq H_T^c=H$, both $R_{q'}$ and $(T^{bc})_{q'}$ lie in $H$. As $K_1$ is a $\pi(n)$-group and $q$ is the least prime divisor of $n$, we have $\pi(|H|)\subseteq\{r\in\pi(n):r>q\}$. Lemma~\ref{grp: orbitbound} now gives $|x^H|\leq s$ for every $x\in\Omega$.

As subgroups of the regular abelian groups $R$ and $T^{bc}$, respectively, $R_{q'}$ and $(T^{bc})_{q'}$ are semiregular. Both have order $s$, since $|R|=|T^{bc}|=q^as$ and $q\nmid s$. Thus every orbit of either subgroup has length $s$. For each $x\in\Omega$, both $x^{R_{q'}}$ and $x^{(T^{bc})_{q'}}$ are contained in $x^H$, whose size is at most $s$. Hence $x^{R_{q'}}=x^H=x^{(T^{bc})_{q'}}$. Take $a_0=bc$ and let $\calE$ be the orbit partition of $H$.

Since $R_{q'}\operatorname{char}R$ and $(T^{a_0})_{q'}\operatorname{char}T^{a_0}$, both $R$ and $T^{a_0}$ permute the blocks of $\calE$, proving~(i). On every $B\in\calE$, the two Hall $q'$-subgroups are transitive and semiregular, so their induced groups are regular. This proves~(ii).

For~(iii), put $K_2=\langle R,T^{a_0}\rangle\leq K_1$. The given sequence is also a complete sequence of block systems for $K_2$. Fix $B\in\calE$ and choose $\omega_0\in B$. Since $R$ is regular, $r\mapsto r(\omega_0)$ is a bijection from $R$ to $\Omega$. Under this identification, $B$ corresponds to $M=R_{q'}$. By Lemma~\ref{grp: normalblocks}, each partition in the sequence is the coset partition of a subgroup of $R$. Thus the sequence corresponds to a subgroup chain $1=S_0<S_1<\cdots<S_\ell=R$ with prime indices $[S_i:S_{i-1}]=p_i$. If a coset of $S_i$ meets $M$ in a point $m$, then its intersection with $M$ is $mS_i\cap M=m(S_i\cap M)$. The restriction of $\calB_i$ to $B$ is therefore the partition of $M$ into cosets of $S_i\cap M$.

Since $R$ is abelian, Lemma~\ref{std: abelian-sylow} gives $S_i=(S_i\cap R_q)\times(S_i\cap M)$. Hence
\[
p_i=[S_i\cap R_q:S_{i-1}\cap R_q]\,[S_i\cap M:S_{i-1}\cap M].
\]
The first factor is a power of $q$ and the second is coprime to $q$. Since $p_i$ is prime,
\[
[S_i\cap M:S_{i-1}\cap M]=
\begin{cases}
1,&p_i=q,\\
p_i,&p_i\ne q.
\end{cases}
\]
Removing repeated terms gives a subgroup chain from $1$ to $M$ with prime indices. The corresponding partitions form a complete sequence of block systems on $B$. When $M=1$, this is the trivial sequence on a singleton.

Finally, each element of $\Stab_{K_2}(B)$ preserves $B$ and each partition $\calB_i$, so it permutes the nonempty intersections $D\cap B$ with $D\in\calB_i$. Consequently, the induced group $\Stab_{K_2}(B)^B$ preserves the restricted sequence, proving~(iii). In particular, its subgroup $\langle R_{q'}^B,(T^{a_0})_{q'}^B\rangle$ also preserves this sequence on $B$.
\end{proof}

\section{Proof of Theorem~\ref{1.1}}\label{sec: fourier}

This section proves Theorem~\ref{1.1}. Before the proof, we record some necessary lemmas.

\begin{lemma}\label{grp: affine}
Let $A$ and $H$ be finite abelian groups of coprime orders, and let $X$ and $Y$ be colored binary Cayley structures on $A\times H$.
Suppose that $f(a,h)=(\beta_h(a),\tau(h))$ is a color-preserving isomorphism from $X$ to $Y$, where each $\beta_h$ is a permutation of $A$ and $\tau\in N_{\Sym(H)}(H_L)$.
Then, for every $h_0\in H$, the map
\[
 f_{h_0}(a,h)=(\beta_{h_0}(a),\tau(h))
\]
is a color-preserving isomorphism from $X$ to $Y$.
\end{lemma}

\begin{proof}
If either factor is trivial, then $f_{h_0}=f$. Assume that both are nontrivial and use additive notation for the group operations.
By Lemma~\ref{std: normalizer}, $\tau(h)=h_*+\theta(h)$ for some $h_*\in H$ and $\theta\in\Aut(H)$.
Put $\eta=\id_A\times\tau$. Since $\eta L_{(a,h)}\eta^{-1}=L_{(a,\theta(h))}$, the relabeled structure $\eta^{-1}Y$ is again a Cayley structure on $A\times H$.
Moreover, $\eta^{-1}f(a,h)=(\beta_h(a),h)$. Replacing $Y$ by $\eta^{-1}Y$ and $f$ by $\eta^{-1}f$, we may therefore assume that $\tau=\id_H$.

For an abelian group $L$, let $(e_x)_{x\in L}$ be the standard basis of $\ell^2(L)$, and let $P_se_x=e_{x+s}$. The adjacency matrix associated with a connection set $S\subseteq L$, whose rows and columns are indexed by $L$, is
\[
 A_S=(\mathbf{1}_{y-x\in S})_{x,y\in L}=\sum_{s\in S}P_{-s}.
\]
For $c\in\calC$, let $D_c^X,D_c^Y\subseteq A\times H$ be the color-$c$ connection sets of $X$ and $Y$, so that
\[
 (x,y)\in R_c^X\iff y-x\in D_c^X,\qquad
 (x,y)\in R_c^Y\iff y-x\in D_c^Y.
\]
For $d\in H$, define the slices
\[
 S_{c,d}=\{a\in A:(a,d)\in D_c^X\},\qquad
 T_{c,d}=\{a\in A:(a,d)\in D_c^Y\}.
\]
Fix a color $c$, and take $a,b\in A$ and $h,d\in H$. Since $f$ is color-preserving and $f(a,h)=(\beta_h(a),h)$, the ordered pair $((a,h),(b,h+d))$ lies in $R_c^X$ exactly when $((\beta_h(a),h),(\beta_{h+d}(b),h+d))$ lies in $R_c^Y$. By the definitions of $D_c^X,D_c^Y$ and the slices, this is equivalent to
\[
 b-a\in S_{c,d}\quad\Longleftrightarrow\quad
 \beta_{h+d}(b)-\beta_h(a)\in T_{c,d}.
\]
Equivalently, the $((a,h),(b,h+d))$-entry of $P_f^{-1}A_{T_{c,d}\times\{d\}}P_f$ is ${\bf 1}_{\beta_{h+d}(b)-\beta_h(a)\in T_{c,d}}$, while the corresponding entry of $A_{S_{c,d}\times\{d\}}$ is ${\bf 1}_{b-a\in S_{c,d}}$; the equivalence above makes these entries equal. Thus $P_f^{-1}A_{T_{c,d}\times\{d\}}P_f=A_{S_{c,d}\times\{d\}}$, that is,
\begin{equation}\label{2}
 P_f A_{S_{c,d}\times\{d\}}P_f^{-1}=A_{T_{c,d}\times\{d\}}.
\end{equation}
Let $\ell$ be a prime not dividing $\exp(A)$. Multiplication by $\ell$ is a permutation of $A$ by Lemma~\ref{std: coprime-map}, so the elements $\ell s$, $s\in S_{c,d}$, are distinct.
The translation matrices commute. Thus, over $\Z/\ell\Z$, we have
\begin{equation}\label{3}
\begin{split}
 A_{S_{c,d}\times\{d\}}^\ell
 &=\left(\sum_{s\in S_{c,d}}P_{-(s,d)}\right)^\ell\\
 &=\sum_{s\in S_{c,d}}P_{-(\ell s,\ell d)}
 =A_{\ell S_{c,d}\times\{\ell d\}}.
\end{split}
\end{equation}
The same calculation gives $A_{T_{c,d}\times\{d\}}^\ell=A_{\ell T_{c,d}\times\{\ell d\}}$ over $\Z/\ell\Z$; here the elements of $\ell T_{c,d}$ are distinct by the same argument.
Taking $\ell$-th powers in~\eqref{2} and applying~\eqref{3} gives a congruence between two $(0,1)$-matrices. Since their entries are $0$ or $1$ and $\ell\ge2$, the two matrices agree over $\Z$, giving
\begin{equation}\label{4}
 P_f A_{\ell S_{c,d}\times\{\ell d\}}P_f^{-1}
 =A_{\ell T_{c,d}\times\{\ell d\}}.
\end{equation}
Since $|H|$ is coprime to $\exp(A)$, we may take $t\ge1$ such that $k=|H|^t\equiv1\pmod{\exp(A)}$.
Let $k=p_1\cdots p_r$, where each $p_j$ is prime and hence coprime to $\exp(A)$. We claim that for every $j=0,\dots,r$,
\[
 P_f A_{m_jS_{c,d}\times\{m_jd\}}P_f^{-1}
 =A_{m_jT_{c,d}\times\{m_jd\}},
 \qquad m_j=p_1\cdots p_j,
\]
with $m_0=1$. For $j=0$ this is~\eqref{2}. Suppose the claim holds for $j-1$. Taking $p_j$-th powers and applying the same calculation as in~\eqref{3} over $\Z/p_j\Z$ gives the claim for $j$; since both sides are $(0,1)$-matrices and $p_j\ge2$, equality over $\Z/p_j\Z$ implies equality over $\Z$. For $j=r$, $m_r=k$, and the claim reads
\begin{equation}\label{5}
 P_f A_{kS_{c,d}\times\{kd\}}P_f^{-1}
 =A_{kT_{c,d}\times\{kd\}}.
\end{equation}
Multiplication by $k$ fixes every element of $A$ and is the zero map on $H$. Thus~\eqref{5} gives
\[
 P_f A_{S_{c,d}\times\{0\}}P_f^{-1}=A_{T_{c,d}\times\{0\}}.
\]
Restricting to the fiber $A\times\{h_0\}$, we find that
\[
 a'-a\in S_{c,d}\quad\Longleftrightarrow\quad
 \beta_{h_0}(a')-\beta_{h_0}(a)\in T_{c,d}
\]
for every $c,d$ and all $a,a'\in A$.
For vertices with second coordinates $h$ and $h'$, take $d=h'-h$. The equivalence shows that $\beta_{h_0}\times\id_H$ preserves every color. Composing on the left with $\eta$ gives the required isomorphism for the original $Y$.
\end{proof}

Now let $A$ and $H$ be finite abelian groups, let $X$ and $Y$ be colored binary Cayley structures on $A\times H$, and let $f(a,h)=(\beta_h(a),\tau(h))$ be a color-preserving isomorphism from $X$ to $Y$ that preserves the $A$-fiber partition.
Index the rows and columns of the adjacency matrix of each color-$c$ relation by the $A$-fibers. For $c\in\calC$ and $d\in H$, the block with row fiber $A\times\{h\}$ and column fiber $A\times\{h+d\}$ is the $|A|\times|A|$ matrix
\begin{equation}\label{four: block-entry}
 A_{c,d}(a,b)=\mathbf 1_{\{((a,h),(b,h+d))\in R_c^X\}}
             =\mathbf 1_{b-a\in S_{c,d}}\qquad(a,b\in A).
\end{equation}
Translation invariance makes this block independent of $h$. Replacing $X$ by $Y$ gives the corresponding block $A'_{c,d}$.
Put $e=\exp(A)$ and $E=\Q(\zeta_e)$. Let $\widehat A$ be the group of complex characters $\chi\colon A\to\mathbb C^\times$, and let $F=(\chi(a))_{a\in A,\,\chi\in\widehat A}$.
By Lemma~\ref{std: fourier}, $F^{-1}=|A|^{-1}\overline F^{\,T}$, and $F$ simultaneously diagonalizes all $A$-circulant matrices. Define $P_he_a=e_{\beta_h(a)}$ and $U_h=F^{-1}P_hF$.

\begin{lemma}\label{four: criterion}
There are scalars $\lambda_{c,d}(\chi),\mu_{c,d}(\psi)\in E$ such that
\[
 F^{-1}A_{c,d}F=\diag(\lambda_{c,d}(\chi)),\qquad
 F^{-1}A'_{c,d}F=\diag(\mu_{c,d}(\psi)).
\]
For $\psi,\chi\in\widehat A$, put $u_h=(U_h)_{\psi,\chi}$. Then
\begin{equation}\label{four: scalar-equation}
 \mu_{c,\tau(k)-\tau(h)}(\psi)u_k=u_h\lambda_{c,k-h}(\chi)
 \qquad(c\in\calC,\ h,k\in H).
\end{equation}
Moreover, $\beta_{h_0}\times\tau$ is a color-preserving isomorphism from $X$ to $Y$ if and only if
\begin{equation}\label{four: frozen-entry}
 \mu_{c,\tau(k)-\tau(h)}(\psi)u_{h_0}=u_{h_0}\lambda_{c,k-h}(\chi)
\end{equation}
for all $\psi,\chi\in\widehat A$, $c\in\calC$, and $h,k\in H$.
\end{lemma}

\begin{proof}
By~\eqref{four: block-entry}, the $(a,b)$-entry of $A_{c,d}$ depends only on $b-a$. Regarding $\chi$ as the column vector with entries $\chi(a)$, we obtain
\[
 (A_{c,d}\chi)(a)=\sum_{s\in S_{c,d}}\chi(a+s)
     =\chi(a)\sum_{s\in S_{c,d}}\chi(s).
\]
Thus the $\chi$-column of $F$ is an eigenvector with eigenvalue $\lambda_{c,d}(\chi)=\sum_{s\in S_{c,d}}\chi(s)\in E$.
Similarly, $A'_{c,d}(a,b)=\mathbf1_{b-a\in T_{c,d}}$, so $A'_{c,d}$ also depends only on $b-a$ and is $A$-circulant. For $\psi\in\widehat A$ and $a\in A$,
\[
(A'_{c,d}\psi)(a)=\sum_{b\in A}A'_{c,d}(a,b)\psi(b)
=\sum_{t\in T_{c,d}}\psi(a+t)
=\psi(a)\sum_{t\in T_{c,d}}\psi(t).
\]
Thus the $\psi$-column of $F$ is an eigenvector of $A'_{c,d}$ with eigenvalue $\mu_{c,d}(\psi)=\sum_{t\in T_{c,d}}\psi(t)\in E$. This proves the two diagonalization formulas.

Fix $h,k\in H$. Since $f$ preserves colors, the block entries satisfy
\begin{equation}\label{four: element-entry}
 A'_{c,\tau(k)-\tau(h)}(\beta_h(a),\beta_k(b))=A_{c,k-h}(a,b).
\end{equation}
These are precisely the $(\beta_h(a),b)$-entries of the two sides of
\[
 A'_{c,\tau(k)-\tau(h)}P_k=P_hA_{c,k-h}.
\]
Conjugating by $F$ gives
\[
 \diag\bigl(\mu_{c,\tau(k)-\tau(h)}(\psi)\bigr)U_k
     =U_h\diag\bigl(\lambda_{c,k-h}(\chi)\bigr).
\]
The diagonal matrices on the left and right scale rows and columns, respectively. The $(\psi,\chi)$-entry is therefore~\eqref{four: scalar-equation}.

A necessary and sufficient condition for $\beta_{h_0}\times\tau$ to preserve color $c$ is that
\[
 A'_{c,\tau(k)-\tau(h)}P_{h_0}=P_{h_0}A_{c,k-h}
\]
for every $h,k\in H$. After conjugation by $F$, this equality holds entrywise exactly when~\eqref{four: frozen-entry} holds for all $\psi,\chi$. Reversing the change of basis proves the converse.
\end{proof}

\begin{lemma}\label{four: propagation}
Let $H$ be a finite abelian group and let $\tau\in\Sym(H)$. Fix $\psi,\chi\in\widehat A$ and let $\lambda_{c,d}=\lambda_{c,d}(\chi)$ and $\mu_{c,d}=\mu_{c,d}(\psi)$.
Suppose that the scalars $(u_h)_{h\in H}$ satisfy $\mu_{c,\tau(k)-\tau(h)}u_k=u_h\lambda_{c,k-h}$ for every $c,h,k$. Let
\[
 \begin{aligned}
 J &=\langle d\in H:\lambda_{c,d}\ne0\text{ for some }c\rangle,\\
 J'&=\langle d\in H:\mu_{c,d}\ne0\text{ for some }c\rangle.
 \end{aligned}
\]
If $u_{h_0}\ne0$, then $u_h\ne0$ for every $h\in h_0+J$, and $\tau(h_0+J)=\tau(h_0)+J'$. In particular, $|J|=|J'|$.
\end{lemma}

\begin{proof}
Put
\[
 \Delta=\{d:\lambda_{c,d}\ne0\text{ for some }c\},\qquad
 \Delta'=\{d:\mu_{c,d}\ne0\text{ for some }c\},
\]
and define digraphs $\Gamma$ and $\Gamma'$ on $H$ by joining $h$ to $h+d$ for each $d\in\Delta$ in $\Gamma$ and each $d\in\Delta'$ in $\Gamma'$. Since $-d=(|H|-1)d$ for every $d\in H$, the additive submonoid generated by $\Delta$ is the subgroup $J$. Thus the vertices reachable from $h$ in $\Gamma$ form the coset $h+J$, and these cosets are precisely the strongly connected components of $\Gamma$. The same argument gives the cosets of $J'$ for $\Gamma'$.

Put $C=h_0+J$. For $h\in C$, take a directed path
\[
 h_0=h^{(0)}\longrightarrow h^{(1)}\longrightarrow\cdots
       \longrightarrow h^{(m)}=h
\]
in $\Gamma$. We prove by induction on $m$ that $u_{h^{(m)}}\ne0$ and $\tau(h^{(m)})\in\tau(h_0)+J'$. For $m=0$, this follows from $u_{h_0}\ne0$ and $\tau(h_0)\in\tau(h_0)+J'$. Suppose the claim holds for $m$. Since $h^{(m)}\to h^{(m+1)}$ is an arc of $\Gamma$, there is $c$ with $\lambda_{c,h^{(m+1)}-h^{(m)}}\ne0$. Taking $h=h^{(m)}$ and $k=h^{(m+1)}$ in~\eqref{four: scalar-equation} gives
\[
\mu_{c,\tau(h^{(m+1)})-\tau(h^{(m)})}u_{h^{(m+1)}}=
u_{h^{(m)}}\lambda_{c,h^{(m+1)}-h^{(m)}}.
\]
The right-hand side is nonzero, so $u_{h^{(m+1)}}\ne0$. Since the coefficient $\mu$ is nonzero, $\tau(h^{(m)})\to\tau(h^{(m+1)})$ is an arc of $\Gamma'$, hence $\tau(h^{(m+1)})-\tau(h^{(m)})\in\Delta'\subseteq J'$. Therefore $\tau(h^{(m+1)})\in\tau(h_0)+J'$, completing the induction. Applying the claim to the chosen path gives $u_h\ne0$ and $\tau(h)\in\tau(h_0)+J'$.

For the reverse inclusion, take $z\in\tau(h_0)+J'$ and a directed path
\[
 z=z^{(0)}\longrightarrow z^{(1)}\longrightarrow\cdots
   \longrightarrow z^{(m)}=\tau(h_0)
\]
in $\Gamma'$. Let $x_i=\tau^{-1}(z^{(i)})$. We prove by backward induction on $i$ that $x_i\in C$ and $u_{x_i}\ne0$. For $i=m$, $x_m=h_0\in C$ and $u_{x_m}=u_{h_0}\ne0$. Suppose the claim holds for $i+1$. Since $z^{(i)}\to z^{(i+1)}$ is an arc of $\Gamma'$, there is $c$ with $\mu_{c,z^{(i+1)}-z^{(i)}}\ne0$. Taking $h=x_i$ and $k=x_{i+1}$ in~\eqref{four: scalar-equation} gives
\[
 u_{x_i}\lambda_{c,x_{i+1}-x_i}
   =\mu_{c,z^{(i+1)}-z^{(i)}}u_{x_{i+1}}.
\]
The left-hand side is nonzero, so $u_{x_i}\ne0$ and $\lambda_{c,x_{i+1}-x_i}\ne0$, hence $x_{i+1}-x_i\in J$. Since $x_{i+1}\in h_0+J$, we get $x_i\in h_0+J=C$. This completes the backward induction. For $i=0$, $x_0\in C$, so $z=\tau(x_0)\in\tau(C)$. We conclude that $\tau(C)=\tau(h_0)+J'$. Comparing cardinalities gives $|J|=|J'|$.
\end{proof}

Let $q$ be an odd prime, let $m$ be a positive integer coprime to $q$, and put $E=\Q(\zeta_m)$. Let $\Lambda$ be a set of size $q^a$, where $a\in\{1,2\}$, and let $Q_0\le\Sym(\Lambda)$ be a regular $q$-group.
Suppose that $Q_0$ and $sQ_0s^{-1}$ lie in a common $q$-subgroup $W_0$ of $\Sym(\Lambda)$, where $s\in\Sym(\Lambda)$. For $t\in\Sym(\Lambda)$, let $P_t$ be the permutation matrix given by $P_te_x=e_{t(x)}$.

\begin{lemma}\label{four: rigidity}
Let $(B_c)_{c\in\calC}$ and $(B'_c)_{c\in\calC}$ be families of matrices over $E$, each of whose matrices commutes with every $P_t$ for $t\in Q_0$. Suppose that the digraph associated with $(B_c)_{c\in\calC}$ is strongly connected.
If an invertible diagonal matrix $D$ over $E$ satisfies $B'_c=P_sDB_cD^{-1}P_s^{-1}$ for every $c\in\calC$, then $D$ is a scalar matrix.
\end{lemma}

\begin{proof}
The field $E$ contains no nontrivial root of unity of $q$-power order. Indeed, if $\zeta_q\in E$, then Lemma~\ref{std: cyclotomic} gives $E=\Q(\zeta_m,\zeta_q)=\Q(\zeta_{mq})$, so
\[
 \varphi(m)=[E:\Q]=\varphi(mq)=(q-1)\varphi(m),
\]
contrary to the assumption that $q$ is odd. Any root of unity of nontrivial $q$-power order has a power of order $q$. It follows that $E^\times$ has no nontrivial finite $q$-subgroup.

Put $M=P_sD$ and consider the monomial matrix group
\[
 \mathcal J=\langle P_t,MP_tM^{-1}:t\in Q_0\rangle.
\]
Replacing the unique nonzero entry in each column by $1$ gives the underlying permutation of a monomial matrix. This defines a homomorphism $\pi\colon\mathcal J\to\Sym(\Lambda)$ with image $\langle Q_0,sQ_0s^{-1}\rangle\le W_0$, a finite $q$-group.
Since $P_t$ commutes with $B_c$, its conjugate $MP_tM^{-1}$ commutes with $B'_c=MB_cM^{-1}$. The matrix $P_t$ also commutes with $B'_c$. Hence every element of $\mathcal J$ commutes with every $B'_c$.

Take $Z=\diag(z_x)\in\ker\pi$. Conjugation by $D$ preserves the support of each $B_c$, and conjugation by $P_s$ permutes both indices according to $s$. Thus $s$ gives an isomorphism of the associated digraphs, and the digraph for $(B'_c)$ is strongly connected.
From $ZB'_c=B'_cZ$ we obtain $z_x=z_y$ whenever $(B'_c)_{x,y}\ne0$. Since the associated digraph of $(B'_c)$ is strongly connected, any two vertices $x,y$ are joined by a directed path $x=x_0\to x_1\to\cdots\to x_m=y$. Applying $z_{x_i}=z_{x_{i+1}}$ successively gives $z_x=z_y$. Hence $Z=zI$. Every $t\in Q_0$ has odd order, so all its cycles have odd length and $\det P_t=1$. Since similarity transformations preserve determinants, $\det(MP_tM^{-1})=\det P_t=1$, and conjugation by $P_s$ also preserves the determinant. Hence every element of $\mathcal J$ has determinant $1$. In particular, $z^{q^a}=\det Z=1$. The preceding argument on roots of unity gives $z=1$. Thus $\pi$ is injective and $\mathcal J$ is a finite $q$-group.

For $x\in\Lambda$, let $\mathcal J_x=\{j\in\mathcal J:\pi(j)(x)=x\}$. Its action on the one-dimensional subspace $Ee_x$ gives a homomorphism $j\mapsto j_{x,x}$ from $\mathcal J_x$ to $E^\times$. The image is a finite $q$-group and is therefore trivial.
For any $j\in\mathcal J$, put $y=\pi(j)(x)$. By regularity there is a unique $t\in Q_0$ with $t(x)=y$. Then $P_t^{-1}j\in\mathcal J_x$, giving $j_{y,x}=(P_t^{-1}j)_{x,x}=1$.
As this holds in every column, every element of $\mathcal J$ is a permutation matrix.

In particular, $DP_tD^{-1}=P_s^{-1}(MP_tM^{-1})P_s$ is a permutation matrix for every $t\in Q_0$. Let $D=\diag(d_x)$. The nonzero entry in column $x$ of $DP_tD^{-1}$ is $d_{t(x)}/d_x$, so $d_{t(x)}=d_x$ for all $x$ and $t$. Transitivity of $Q_0$ now gives the result.
\end{proof}

\begin{lemma}\label{four: cyclic-normalizer}
Let $q$ be a prime and let $H=C_{q^2}$. Let each $h\in H$ be written as $i+qj$, with $0\leq i,j<q$, and identify $h$ with the pair $(i,j)\in\F_q^2$. To add $(i,j)$ and $(i',j')$, take $(i+i',j+j')$ if $i+i'<q$ and $(i+i'-q,j+j'+1)$ otherwise, reducing the second coordinate modulo $q$. Let $Q=H_L$ and put
\begin{equation}
 W=\{(i,j)\mapsto(i+t,j+b_i): t\in\F_q,\ (b_i)\in\F_q^q\}.
 \label{6}
\end{equation}
Then $W$ is the unique Sylow $q$-subgroup of $\Sym(H)$ containing $Q$, and $N_{\Sym(H)}(W)$ consists precisely of the permutations of the form
\begin{equation}
 \tau(i,j)=(ci+d,aj+b_i),\qquad
 a,c\in\F_q^\times,\quad d\in\F_q,\quad(b_i)\in\F_q^q.
 \label{four: cyclic-normalizer-form}
\end{equation}
Moreover, every $\tau\in N_{\Sym(H)}(W)$ satisfies
\begin{equation}
 \tau(h+qt)-\tau(h)=qat\qquad(h\in H,\ t\in\F_q),
 \label{four: uniform-slope}
\end{equation}
where $a$ is as in~\eqref{four: cyclic-normalizer-form}. In particular, if $Q,\tau Q\tau^{-1}\leq W$, then $\tau\in N_{\Sym(H)}(W)$.
\end{lemma}
\begin{proof}
By Lemma~\ref{std: legendre}, a Sylow $q$-subgroup of $\Sym(H)$ has order $q^{q+1}$. The permutations in~\eqref{6} form a group of this order, so $W$ is a Sylow $q$-subgroup. The translation $h\mapsto h+1$, which generates $Q$, acts by $(i,j)\mapsto(i+1,j+\delta_{i,q-1})$. Here the first coordinate is taken modulo $q$, and the Kronecker delta records the carry into the second coordinate. This translation belongs to $W$, so $Q\leq W$.

Let $V$ be any Sylow $q$-subgroup of $\Sym(H)$ containing $Q$. Since $V$ is a nontrivial $q$-group, its center contains a subgroup $Z$ of order $q$. By Lemma~\ref{std: centralizer}, $Z\leq Q$. The cyclic group $Q$ has a unique subgroup of order $q$, so $Z$ is generated by $h\mapsto h+q$. Its orbits are $B_i=\{(i,j):j\in\F_q\}$, for $i\in\F_q$. Since $Z\leq Z(V)$, the group $V$ preserves the partition $D=\{B_i:i\in\F_q\}$, and every $v\in V$ commutes with $(i,j)\mapsto(i,j+1)$. Hence $v(i,j)=(\sigma(i),j+b_i)$ for some $\sigma\in\Sym(\F_q)$ and $(b_i)\in\F_q^q$.
The group $V^D$ induced on the blocks is a $q$-group containing the cyclic group generated by $i\mapsto i+1$. By Lemma~\ref{std: legendre}, $v_q(q!)=1$, so $|V^D|\leq q$. The cyclic subgroup already has order $q$ and therefore equals $V^D$. Thus $V\leq W$, and equality follows by comparing orders.

Every $W$-invariant partition into blocks of size $q$ is also $Q$-invariant. The regular cyclic group $Q$ has just one such partition. Hence every $\tau\in N_{\Sym(H)}(W)$ preserves $D$. Its induced action on $D$ normalizes the cyclic group generated by $i\mapsto i+1$ and so has the form $i\mapsto ci+d$. The subgroup $\operatorname{fix}_W(D)$ fixing each $B_i$ setwise induces on $B_i$ the cyclic group generated by $(i,j)\mapsto(i,j+1)$. Conjugation by $\tau$ permutes these induced groups, so
\[
 \tau(i,j)=(ci+d,a_i j+b_i),\qquad a_i\in\F_q^\times.
\]
Let $\sigma_1(i,j)=(i+1,j)$, which lies in $W$. Since $\tau\sigma_1\tau^{-1}\in W$, its action on the second coordinate is a translation on each block. Its value at $\tau(i,j)$ is $\tau(i+1,j)=(c(i+1)+d,a_{i+1}j+b_{i+1})$. Comparing coefficients of $j$ gives $a_{i+1}=a_i$. Thus all the $a_i$ have the same nonzero value $a$, as claimed in~\eqref{four: cyclic-normalizer-form}.

Conversely, suppose that $\tau$ has the form~\eqref{four: cyclic-normalizer-form}. With $x=ci+d$ and $y=aj+b_i$, we have
\[
 \tau\sigma_1\tau^{-1}(x,y)=(x+c,y+b_{i+1}-b_i),
 \qquad i=c^{-1}(x-d).
\]
This conjugate belongs to $W$. Conjugation by $\tau$ takes a translation in the second coordinate supported on $B_i$ to a translation supported on $B_{ci+d}$ and multiplies the translation parameter by $a$. These blockwise translations and $\sigma_1$ generate $W$. Thus $\tau W\tau^{-1}\leq W$, with equality by comparison of orders. Finally, $h$ and $h+qt$ have the same first coordinate, and their second coordinates differ by $t$. The second coordinates of their images differ by $at$, which gives~\eqref{four: uniform-slope} in $C_{q^2}$.

If $Q,\tau Q\tau^{-1}\leq W$, then $W$ and $\tau W\tau^{-1}$ are both Sylow $q$-subgroups containing the regular cyclic group $\tau Q\tau^{-1}$. The uniqueness just proved, applied to this group, gives $\tau W\tau^{-1}=W$.
\end{proof}

\begin{lemma}\label{four: shear}
Let $q$ be a prime, let $H=C_q^2$, and let $Q=H_L$. If $W$ is a Sylow $q$-subgroup of $\Sym(H)$ containing $Q$, then there are linear coordinates on $H=\F_q^2$ in which
\begin{equation}
 W=\{(i,j)\mapsto(i+t,j+b_i): t\in\F_q,\ (b_i)\in\F_q^q\}.
 \label{7}
\end{equation}
Every regular elementary abelian subgroup $T\leq W$ has the form $T=wQw^{-1}$ for some permutation
\begin{equation}
 w(i,j)=(i,j+c_i),\qquad(c_i)\in\F_q^q.
 \label{four: shear-form}
\end{equation}
In particular, if $\tau Q\tau^{-1}\leq W$, then
\begin{equation}
 \tau=wn,\qquad n\in N_{\Sym(H)}(Q)=\AGL(2,q).
 \label{four: shear-factorization}
\end{equation}
\end{lemma}
\begin{proof}
Choose a subgroup $Z_0\leq Z(W)$ of order $q$. Since $Q\leq W$, the group $Z_0$ centralizes $Q$. By Lemma~\ref{std: centralizer}, $C_{\Sym(H)}(Q)=Q$. Hence $Z_0\leq Q$, and $Z_0$ consists of the translations by the vectors of a one-dimensional subspace of $H$. Choose linear coordinates in which this subspace is $V=\{(0,j):j\in\F_q\}$. Then $Z_0$ is generated by $(i,j)\mapsto(i,j+1)$.

Every $w_0\in W$ commutes with $Z_0$, so it permutes the $Z_0$-orbits $B_i=\{(i,j):j\in\F_q\}$ and has the form $w_0(i,j)=(\sigma(i),j+b_i)$. The group $W^D$ induced on $D=\{B_i:i\in\F_q\}$ is a $q$-group containing the group generated by $i\mapsto i+1$. Since $v_q(q!)=1$ by Lemma~\ref{std: legendre}, $W^D$ has order at most $q$ and hence equals this cyclic group. Thus $W$ is contained in the group on the right-hand side of~\eqref{7}. Both groups have order $q^{q+1}$, so they coincide.

Let $T\leq W$ be regular and elementary abelian, and let $\pi\colon W\to\Sym(D)$ be the permutation representation of $W$ on $D$. Since $T$ is transitive, its image $\pi(T)$ is transitive on $D$ and thus equals the group of order $q$ generated by $i\mapsto i+1$. It follows that $Z_T=\ker(\pi|_T)$ has order $q$. Choose $t\in T$ inducing $i\mapsto i+1$ on $D$. Then $t(i,j)=(i+1,j+b_i)$.
By~\eqref{7}, each $z\in Z_T$ has the form $z(i,j)=(i,j+z_i)$. Comparing second coordinates in $zt=tz$ gives $z_{i+1}=z_i$. Hence $Z_T$ is contained in $\{(i,j)\mapsto(i,j+z):z\in\F_q\}$, and both groups have order $q$. Moreover, $T=\langle t,Z_T\rangle$ and $t^q=1$ because $T$ is elementary abelian. Since $t^q(i,j)=(i,j+\sum_{r=0}^{q-1}b_r)$, we have
\begin{equation}
 \sum_{i=0}^{q-1}b_i=0.
 \label{8}
\end{equation}
Let $c_0=0$ and $c_i=\sum_{r=0}^{i-1}b_r$ for $1\leq i<q$. By~\eqref{8}, the identity
\begin{equation}
 c_{i+1}-c_i=b_i\qquad(i\in\F_q)
 \label{four: shear-difference}
\end{equation}
also holds for $i=q-1$, with indices read modulo $q$. Define $w(i,j)=(i,j+c_i)$. Then $w^{-1}tw(i,j)=(i+1,j+c_i+b_i-c_{i+1})=(i+1,j)$, and $w$ commutes with every element of $Z_T$. Thus $w^{-1}Tw$ is generated by $(i,j)\mapsto(i+1,j)$ and $Z_T$. These generate $Q$, proving $T=wQw^{-1}$.

Now apply this result to $T=\tau Q\tau^{-1}\leq W$. We have $w^{-1}\tau Q\tau^{-1}w=Q$, so $n=w^{-1}\tau\in N_{\Sym(H)}(Q)=\AGL(2,q)$ by Lemma~\ref{std: normalizer}. This gives~\eqref{four: shear-factorization}.
\end{proof}

\begin{proof}[Proof of Theorem~\ref{1.1}]
Since $f$ preserves the $A$-fiber partition, $f(\calE_A)=\calE_A$. Both $R$ and $fRf^{-1}$ preserve $f(\calE_A)$. Denote the group induced by $B=\langle R,fRf^{-1}\rangle$ on this partition by $\overline B$, and put $Q=H_L$. Identify the fiber $A\times\{h\}$ with $h$. Then $R$ induces $Q$ on the fibers. Since $f(a,h)=(\beta_h(a),\tau(h))$, the conjugate $fRf^{-1}$ induces $\tau Q\tau^{-1}$. Thus $\overline B$ contains both groups. By Lemma~\ref{std: sylow}(ii), some Sylow $q$-subgroup $S$ of $\overline B$ contains $Q$. By Lemma~\ref{std: sylow}(iii), there is an element $\overline b\in\overline B$ such that $\overline b\tau Q\tau^{-1}\overline b^{-1}\leq S$. Choose a preimage $b$ of $\overline b$ in $B$. Since $b\in\Aut(Y)$, the map $bf$ is again a color-preserving isomorphism from $X$ to $Y$ that preserves the fiber partition. Retain the notation $f$ for $bf$ and let $f(a,h)=(\beta_h(a),\tau(h))$ denote its coordinate maps. By Lemma~\ref{std: sylow}(ii), there is a Sylow $q$-subgroup $W$ of $\Sym(H)$ containing $S$. Thus
\begin{equation}
 Q,\tau Q\tau^{-1}\leq W.
 \label{four: common-sylow}
\end{equation}

If $|H|=q$, then $W$ has order $q$ and hence equals $Q$. If $H=C_4$, then $W$ is the dihedral group of order eight and $Q$ is its unique cyclic subgroup of order four. If $H=C_2^2$, then $Q$ is the normal Klein four-group $V_4=\{1,(12)(34),(13)(24),(14)(23)\}$ in $\Sym(H)\cong\Sym(4)$. In each case, $\tau Q\tau^{-1}=Q$, and Lemma~\ref{grp: affine} shows that $\beta_{h_0}\times\tau$ is a color-preserving isomorphism for every $h_0\in H$.

Now suppose that $q$ is odd and $|H|=q^2$. Put $E=\Q(\zeta_{\exp(A)})$. Since $\gcd(|A|,|H|)=1$ and $\exp(A)\mid |A|$, we have $q\nmid\exp(A)$. Fix $h_0\in H$. For each pair of characters $(\psi,\chi)\in\widehat A\times\widehat A$, let $u_h=(U_h)_{\psi,\chi}$, $\lambda_{c,d}=\lambda_{c,d}(\chi)$ and $\mu_{c,d}=\mu_{c,d}(\psi)$. Define $J$ and $J'$ as in Lemma~\ref{four: propagation}. Suppose first that $u_{h_0}\ne0$ and consider~\eqref{four: scalar-equation} for this pair.

If $J=H$, Lemma~\ref{four: propagation} gives $J'=H$ and $u_h\ne0$ for every $h\in H$. Put $B_c(h,k)=\lambda_{c,k-h}$, $B'_c(h,k)=\mu_{c,k-h}$ and $D_u=\diag(u_h)$. Both matrix families $(B_c)_{c\in\calC}$ and $(B'_c)_{c\in\calC}$ commute with the permutation matrices of $Q$, and the digraph associated with $(B_c)_{c\in\calC}$ is strongly connected. In~\eqref{four: scalar-equation}, replace $h,k$ by $\tau^{-1}(h),\tau^{-1}(k)$ and divide by $u_{\tau^{-1}(k)}$. This gives
\[
 \mu_{c,k-h}
 =u_{\tau^{-1}(h)}\lambda_{c,\tau^{-1}(k)-\tau^{-1}(h)}
 u_{\tau^{-1}(k)}^{-1}.
\]
The expression on the right is the $(\tau^{-1}(h),\tau^{-1}(k))$-entry of $D_uB_cD_u^{-1}$ and hence the $(h,k)$-entry of $P_\tau D_uB_cD_u^{-1}P_\tau^{-1}$. Thus $B'_c=P_\tau D_uB_cD_u^{-1}P_\tau^{-1}$. By~\eqref{four: common-sylow} and Lemma~\ref{four: rigidity}, $D_u$ is a scalar matrix, and therefore
\begin{equation}
 \mu_{c,\tau(k)-\tau(h)}=\lambda_{c,k-h}
 \qquad(c\in\calC,\ h,k\in H).
 \label{four: equality-all-pairs}
\end{equation}
If $J=\{0\}$, then $J'=\{0\}$ by Lemma~\ref{four: propagation}. Taking $h=k=h_0$ in~\eqref{four: scalar-equation} gives $\mu_{c,0}=\lambda_{c,0}$. For $k\ne h$, the two coefficients in~\eqref{four: equality-all-pairs} vanish since $J=J'=\{0\}$. Thus~\eqref{four: equality-all-pairs} holds in this case as well.

Suppose that $H=C_{q^2}$, and use additive notation for its group operation, so the identity is $0$ and the $t$-th power of $h$ is $th$. Its unique subgroup of order $q$ is $qH=\{qt:t\in\F_q\}$, generated by the element $q=q\cdot1$; its cosets are $q$ blocks of size $q$. By~\eqref{four: common-sylow} and Lemma~\ref{four: cyclic-normalizer}, $\tau\in N_{\Sym(H)}(W)$. Its restriction to each coset of $qH$ has the same multiplier $a\in\F_q^\times$. In particular, $\tau(h_0+qt)=\tau(h_0)+qat$ for every $t\in\F_q$. It remains to consider $|J|=q$, in which case Lemma~\ref{four: propagation} gives $J=J'=qH$. Parametrize $h_0+qH$ by $h_0+qt$ and its image by $\tau(h_0)+qt'$. In these coordinates the restriction of $\tau$ is $t\mapsto at$. Identifying both cosets with $\F_q$, let
\[
 C_c(t,r)=\lambda_{c,q(r-t)},\quad
 C'_c(t,r)=\mu_{c,q(r-t)},\quad
 D_0=\diag_{t\in\F_q}(u_{h_0+qt}).
\]
By Lemma~\ref{four: propagation}, $D_0$ is invertible. Let $s_a(t)=at$. Equation~\eqref{four: scalar-equation} gives $C'_c=P_{s_a}D_0C_cD_0^{-1}P_{s_a}^{-1}$. Both matrix families commute with the regular translation group $Q_0$ on $\F_q$. Since $J=qH$, the digraph associated with $(C_c)_{c\in\calC}$ is strongly connected. Also, $s_aQ_0s_a^{-1}=Q_0$. Apply Lemma~\ref{four: rigidity} with $W_0=Q_0$ to conclude that $D_0$ is a scalar matrix. Hence $u_h$ is constant on $h_0+qH$, and
\begin{equation}
 \mu_{c,qat}=\lambda_{c,qt}\qquad(c\in\calC,\ t\in\F_q).
 \label{four: cyclic-component-coefficients}
\end{equation}
For $k-h\in qH$, equations~\eqref{four: uniform-slope} and~\eqref{four: cyclic-component-coefficients} give~\eqref{four: equality-all-pairs}. For $k-h\notin qH$, the coset partition is $\tau$-invariant, so $\tau(k)-\tau(h)\notin qH$. Since $J=J'=qH$, both coefficients in~\eqref{four: equality-all-pairs} vanish.

Now let $H=C_q^2$. By Lemma~\ref{four: shear}, we may choose linear coordinates and let $\tau=wn$, where $w(i,j)=(i,j+c_i)$ and $n\in\AGL(2,q)$. Put $\eta=\id_A\times n$ and $X^*=\eta X$. By Lemma~\ref{std: normalizer}, $\AGL(2,q)=N_{\Sym(H)}(H_L)$. Thus $n$ normalizes $Q$, and $\eta$ normalizes $(A\times H)_L$. Consequently, $X^*$ is again a Cayley structure on $A\times H$. The isomorphism $f^*=f\eta^{-1}\colon X^*\to Y$ has the form $f^*(a,z)=(\gamma_z(a),w(z))$, where $\gamma_z=\beta_{n^{-1}(z)}$.
For $X^*,Y,f^*$, define the blocks $A_{c,d},A'_{c,d}$ by~\eqref{four: block-entry}, the coefficients $\lambda,\mu,u_z$ as in Lemma~\ref{four: criterion}, and the subgroups $J,J'$ as in Lemma~\ref{four: propagation}, retaining the same notation in the new coordinates. Put $z_0=n(h_0)$ and $V=\{(0,j):j\in\F_q\}$. Then $\gamma_{z_0}=\beta_{h_0}$. Since $nQn^{-1}=Q$ and $\tau=wn$, we still have $Q,wQw^{-1}\leq W$.

Suppose that some character pair with $u_{z_0}\ne0$ gives a one-dimensional subgroup $J$ distinct from $V$. Let $z_0+J=\{(i,ri+s):i\in\F_q\}$. Its image under $w$ is $\{(i,ri+s+c_i):i\in\F_q\}$. By Lemma~\ref{four: propagation}, this image is the coset $w(z_0)+J'$. Since $|J'|=|J|=q$, the subgroup $J'$ is one-dimensional, and the points in the image have the form $(i,r'i+s')$ as $i$ ranges over $\F_q$. Comparing second coordinates gives $ri+s+c_i=r'i+s'$ for all $i$. Thus $c_i=(r'-r)i+(s'-s)$ and $w$ is affine. Lemma~\ref{grp: affine} now gives an isomorphism $\gamma_{z_0}\times w\colon X^*\to Y$.

Otherwise, every character pair with $u_{z_0}\ne0$ has $J\in\{\{0\},V,H\}$. If $J=H$, the argument leading to~\eqref{four: equality-all-pairs} applies with $w$ in place of $\tau$, since $Q,wQw^{-1}\leq W$. In particular, Lemma~\ref{four: rigidity} makes $\diag(u_z)$ a scalar matrix, giving $\mu_{c,w(k)-w(h)}=\lambda_{c,k-h}$ for all $h,k\in H$. If $J=\{0\}$, then $J'=\{0\}$ by Lemma~\ref{four: propagation}. Thus $\lambda_{c,d}=\mu_{c,d}=0$ for $d\ne0$, and both sides of~\eqref{four: equality-all-pairs} vanish when $h\ne k$. Taking $h=k=z_0$ in~\eqref{four: scalar-equation} gives $\mu_{c,0}u_{z_0}=u_{z_0}\lambda_{c,0}$; since $u_{z_0}\ne0$, we obtain $\mu_{c,0}=\lambda_{c,0}$. These coefficients do not depend on the diagonal position, so~\eqref{four: equality-all-pairs} also holds for $h=k$.

If $J=V$, let $z_0+V=\{(i_0,j):j\in\F_q\}$, where $i_0$ is the first coordinate of $z_0$. Lemma~\ref{four: propagation} gives $w(z_0+V)=w(z_0)+J'$. Since $w(i,j)=(i,j+c_i)$ fixes each coset of $V$ setwise, $w(z_0+V)=w(z_0)+V$ and hence $J'=V$. Parametrize $z_0+V$ and $w(z_0)+V$ by assigning the same parameter $t\in\F_q$ to $z_0+(0,t)$ and $w(z_0)+(0,t)$. In these coordinates the restriction of $w$ is the identity. Put
\[
 C_c(t,r)=\lambda_{c,(0,r-t)},\quad
 C'_c(t,r)=\mu_{c,(0,r-t)},\quad
 D_0=\diag_{t\in\F_q}(u_{z_0+(0,t)}).
\]
By Lemma~\ref{four: propagation}, $D_0$ is invertible, and~\eqref{four: scalar-equation} gives $C'_c=D_0C_cD_0^{-1}$. Both matrix families commute with the regular translation group $Q_0=(\F_q)_L$. Since $J=V$, the digraph associated with $(C_c)_{c\in\calC}$ is strongly connected. Lemma~\ref{four: rigidity}, with $s=\id$ and $W_0=Q_0$, shows that $D_0$ is a scalar matrix. It follows that $\mu_{c,d}=\lambda_{c,d}$ for $c\in\calC$ and $d\in V$.
If $k=(i_1,j_1)$ and $h=(i_2,j_2)$ lie in the same coset of $V$, then $i_1=i_2$ and $w(k)-w(h)=(0,j_1-j_2)=k-h\in V$. If $i_1\ne i_2$, both $k-h$ and $w(k)-w(h)$ lie outside $V$. Since $J=J'=V$, the corresponding coefficients vanish. This proves~\eqref{four: equality-all-pairs} for all $h,k\in H$.

Fix a character pair $(\psi,\chi)$ and a color $c$. If $u_{h_0}\ne0$, multiply~\eqref{four: equality-all-pairs} by $u_{h_0}$ to obtain~\eqref{four: frozen-entry}. If $u_{h_0}=0$, both sides of~\eqref{four: frozen-entry} vanish. Thus~\eqref{four: frozen-entry} holds for every character pair, every color and all $h,k\in H$ (in the elementary abelian case, $z_0$ denotes the image of $h_0$ in the new coordinates). By Lemma~\ref{four: criterion}, $\beta_{h_0}\times\tau\colon X\to Y$ is a color-preserving isomorphism in the cyclic case, as is $\gamma_{z_0}\times w\colon X^*\to Y$ in the elementary abelian case. Composing the latter map with $\eta$ yields $(\gamma_{z_0}\times w)\eta=\beta_{h_0}\times\tau$, since $wn=\tau$ and $\gamma_{z_0}=\beta_{h_0}$. Thus $(a,h)\mapsto(\beta_{h_0}(a),\tau(h))$ is a color-preserving isomorphism from $X$ to $Y$ in this case as well. The choice of $b\in B$ preceded that of $h_0$, so the same $b$ works for every $h_0\in H$.
\end{proof}

\section{Proofs of Theorems~\ref{thm: main} and~\ref{thm: solvingset}}\label{sec: proof}

This section proves Theorems~\ref{thm: main} and~\ref{thm: solvingset}. We first prove Theorem~\ref{thm: main}, then prove Lemma~\ref{grp: coloursubstitution}, and finally deduce Theorem~\ref{thm: solvingset} from these two results.

\begin{proof}[Proof of Theorem~\ref{thm: main}]
We induct on the number $r$ of distinct prime divisors of $|G|$. If $G$ is trivial, then $r=0$ and the identity map on the single vertex, which is common to $X$ and $Y$, is the required color-preserving isomorphism. If $r=1$, then $G$ itself is its only Sylow subgroup. $D_G=\Sym(G)$, so every color-preserving isomorphism from $X$ to $Y$ lies in $D_G$. Suppose that $r\geq2$ and that the theorem holds for abelian groups of cube-free order with fewer than $r$ distinct prime divisors. Let $q=\min\pi(|G|)$ and let $G=A\times H$, where $A=\prod_{p>q}G_p$ and $H=G_q$. Since $|G|$ is cube-free, $H$ is isomorphic to $C_q$, $C_{q^2}$ or $C_q^2$.

Put $R=G_L$. The group $B=\langle R,fRf^{-1}\rangle\leq\Aut(Y)$ is conjugate to $\langle G_L,f^{-1}G_Lf\rangle$. Conjugating a complete sequence of block systems for $\langle G_L,f^{-1}G_Lf\rangle$ therefore gives one for $B$. By Lemma~\ref{grp: alignment}, there exists $a_0\in B$ such that $A_L$ and $(fA_Lf^{-1})^{a_0}$ have the same orbit partition, that is, $a_0^{-1}(f(\calE_A))=\calE_A$. Thus $f_1=a_0^{-1}f$ preserves the $A$-fiber partition and has the form $f_1(a,h)=(\beta_h(a),\tau(h))$. The group $B_1=\langle R,f_1Rf_1^{-1}\rangle\leq B$ preserves every partition in the complete sequence as well as the fiber partition $f_1(\calE_A)$.

By Theorem~\ref{1.1}, there exists $b\in B_1$ such that $f_2=bf_1$ has the following property: writing $f_2(a,h)=(\gamma_h(a),\rho(h))$, for every fixed $h_0\in H$, the choice $\alpha=\gamma_{h_0}$ gives a color-preserving isomorphism
\begin{equation}\label{grp: frozeniso}
 (\alpha,\rho)\colon X\longrightarrow Y.
\end{equation}
Since $B_2=\langle R,f_2Rf_2^{-1}\rangle\leq B_1$, the group $B_2$ also preserves all the partitions in the complete sequence and the fiber partition $f_1(\calE_A)$.

Define colored structures $U,V$ on $A$, with relations indexed by $(c,h,k)\in\calC\times H^2$, by
\begin{align*}
 R^U_{c,h,k}
 &=\{(a,a'): \bigl((a,h),(a',k)\bigr)\in R_c^X\},\\
 R^V_{c,h,k}
 &=\{(a,a'): \bigl((a,\rho(h)),(a',\rho(k))\bigr)\in R_c^Y\}.
\end{align*}
These relations are invariant under the translation group $A_L$, so $U,V$ are colored binary Cayley structures on $A$. By~\eqref{grp: frozeniso}, $\alpha\colon U\to V$ is a color-preserving isomorphism.

Put $C=A\times\{\rho(h_0)\}$. Both $A_L$ and $f_2A_Lf_2^{-1}$ stabilize $C$, since $A_L$ fixes each fiber setwise and $f_2$ maps $A\times\{h_0\}$ onto $C$. Identify $C$ with $A$ via the first coordinate. On this fiber $f_2$ induces $\alpha$, so the two groups induce $A_L$ and $\alpha A_L\alpha^{-1}$ on $C$. In particular,
\[
 \langle A_L,\alpha A_L\alpha^{-1}\rangle
 \leq \Stab_{B_2}(C)^C \leq \Stab_{B_1}(C)^C.
\]
Since $C$ is a block of the fiber partition $f_1(\calE_A)$, Lemma~\ref{grp: alignment}(iii) shows that the group on the right preserves a complete sequence of block systems on $C$. The group on the left preserves the same sequence. Conjugating by $\alpha$ therefore gives a complete sequence for $\langle A_L,\alpha^{-1}A_L\alpha\rangle$.

The order of $A$ is cube-free and has $r-1$ distinct prime divisors. By induction, there is a color-preserving isomorphism $s=\prod_{p>q}s_p\colon U\to V$ with $s_p\in\Sym(G_p)$. By construction, $s$ maps every $R^U_{c,h,k}$ onto $R^V_{c,h,k}$. Thus $(s,\rho)\colon X\to Y$ is the required color-preserving isomorphism, with $(s,\rho)\in\prod_{p\mid |G|}\Sym(G_p)$.
\end{proof}

\begin{lemma}\label{grp: coloursubstitution}
Let $G=G_1\times\cdots\times G_r$ be a direct product of finite groups, and let $X,Y$ be colored binary Cayley structures on $G$. For each $i$, let $S_i$ be a universal color-preserving solving set for $G_i$. Thus, for any finite families $(Q_\alpha)$ and $(Q'_\alpha)$ of $(G_i)_L$-invariant binary relations in $G_i^2$, if a permutation $\sigma$ of $G_i$ maps $Q_\alpha$ onto $Q'_\alpha$ for every $\alpha$, then some $s_i\in S_i$ does so simultaneously for all $\alpha$. Suppose that $h=(h_1,\dots,h_r)\in\prod_i\Sym(G_i)$ is a color-preserving isomorphism from $X$ to $Y$, acting by
\[
 h(x_1,\dots,x_r)=(h_1(x_1),\dots,h_r(x_r)),\qquad h_i\colon G_i\to G_i.
\]
Then there is an element $(s_1,\dots,s_r)\in\prod_i S_i$ such that $s=(s_1,\dots,s_r)\in\prod_i\Sym(G_i)$ is a color-preserving isomorphism from $X$ to $Y$.
\end{lemma}
\begin{proof}
Fix $i$ and put $G_{-i}=\prod_{j\ne i}G_j$. Let elements of $G$ be written as $(a,u)$, with $a$ in the $i$th factor. For every color $c$ and every $(u,v)\in G_{-i}^2$, define a binary relation on $G_i$ by
\begin{equation}\label{colour: source}
 Q_{c,u,v}=\{(a,b)\in G_i^2: ((a,u),(b,v))\in R_c^X\}.
\end{equation}
For a map $\phi=(\phi_1,\dots,\phi_r)$ on $G$, put $\phi_{-i}=\prod_{j\ne i}\phi_j$, and define
\begin{equation}\label{colour: target}
 Q'_{c,u,v}(\phi)=\{(a,b)\in G_i^2:
 ((a,\phi_{-i}(u)),(b,\phi_{-i}(v)))\in R_c^Y\}.
\end{equation}
Since $X,Y$ are Cayley structures on $G$, we have $G_L\leq\Aut(X)$ and $G_L\leq\Aut(Y)$. For any $g\in G_i$, left multiplication by $(1,\dots,g,\dots,1)$ on both entries preserves membership in every $R_c^X$ and $R_c^Y$. It follows that $Q_{c,u,v}$ and $Q'_{c,u,v}(\phi)$ are $(G_i)_L$-invariant.

Let $\phi=\prod_{j=1}^r\phi_j\in\prod_j\Sym(G_j)$ be a color-preserving isomorphism from $X$ to $Y$. For all $a,b\in G_i$ we have
\begin{equation*}\begin{aligned}
 (a,b)\in Q_{c,u,v}
 &\iff ((a,u),(b,v))\in R_c^X\\
 &\iff ((\phi_i(a),\phi_{-i}(u)),(\phi_i(b),\phi_{-i}(v)))\in R_c^Y\\
 &\iff (\phi_i(a),\phi_i(b))\in Q'_{c,u,v}(\phi).
\end{aligned}\end{equation*}
Thus $\phi_i$ maps each relation $Q_{c,u,v}$ defined by~\eqref{colour: source} onto the corresponding relation $Q'_{c,u,v}(\phi)$ defined by~\eqref{colour: target}. By the universality of $S_i$, there is an $s_i\in S_i$ such that $s_i(Q_{c,u,v})=Q'_{c,u,v}(\phi)$ for every $c,u,v$. Hence, for every color-preserving isomorphism $\phi=\prod_{j=1}^r\phi_j\in\prod_j\Sym(G_j)$ from $X$ to $Y$ and every $i$, such an $s_i$ exists.

Let $\phi^{(s)}$ be the map obtained from $\phi$ by replacing its $i$th factor with $s_i$:
\[
 \phi^{(s)}=(\phi_1,\dots,\phi_{i-1},s_i,\phi_{i+1},\dots,\phi_r).
\]
Then for all $a,b\in G_i$ and $u,v\in G_{-i}$,
\begin{equation*}\begin{aligned}
 ((a,u),(b,v))\in R_c^X
 &\iff (a,b)\in Q_{c,u,v}\\
 &\iff (s_i(a),s_i(b))\in Q'_{c,u,v}(\phi)\\
 &\iff ((s_i(a),\phi_{-i}(u)),(s_i(b),\phi_{-i}(v)))\in R_c^Y.
\end{aligned}\end{equation*}
Thus $\phi^{(s)}$ is again a color-preserving isomorphism from $X$ to $Y$.

We now replace the factors successively. For $i=1,\dots,r$, we construct a color-preserving isomorphism $\phi^{(i)}=(\phi^{(i)}_1,\dots,\phi^{(i)}_r)\in\prod_j\Sym(G_j)$ from $X$ to $Y$ such that $\phi^{(i)}_j\in S_j$ for $j\le i$. Let $\phi^{(0)}=h$, so that $\phi^{(0)}\in\prod_j\Sym(G_j)$ and $\phi^{(0)}$ is a color-preserving isomorphism from $X$ to $Y$. Suppose that for some $i$ we have constructed a color-preserving isomorphism $\phi^{(i-1)}=(\phi^{(i-1)}_1,\dots,\phi^{(i-1)}_r)\in\prod_j\Sym(G_j)$ from $X$ to $Y$ with $\phi^{(i-1)}_j\in S_j$ for $j\le i-1$. Apply the preceding paragraph to $\phi=\phi^{(i-1)}$: there is $s_i\in S_i$ such that replacing the $i$th factor of $\phi^{(i-1)}$ by $s_i$ gives a color-preserving isomorphism $\phi^{(i)}$ from $X$ to $Y$, and $\phi^{(i)}\in\prod_j\Sym(G_j)$. Then $\phi^{(i)}_j\in S_j$ for $j\le i$. For $i=r$, we obtain $\phi^{(r)}=(s_1,\dots,s_r)\in\prod_j\Sym(G_j)$ with $s_j=\phi^{(r)}_j\in S_j$ for every $j$, and $\phi^{(r)}\colon X\to Y$ is a color-preserving isomorphism.
\end{proof}

\begin{proof}[Proof of Theorem~\ref{thm: solvingset}]
Let $C=\{c_0\}$, and let $R_{c_0}$ be the arc set $\{(g,gs):g\in G,\ s\in S\}$ of $X=\Cay(G,S)$. Then $X=(G,(R_c)_{c\in C})$ is a colored binary Cayley structure with color set $C$, where $R_{c_0}\subseteq G^2$ is the arc relation of $X$ and $G_L\le\operatorname{Aut}(X)$ because $X$ is a Cayley digraph. For any Cayley digraph $Y$ on $G$ isomorphic to $X$, let $R'_{c_0}$ be its arc set. Then $Y=(G,(R'_c)_{c\in C})$ is a colored binary Cayley structure with the same color set $C$, and a color-preserving isomorphism $X\to Y$ is exactly a digraph isomorphism $X\to Y$. Let $Y$ be any such Cayley digraph.

The hypothesis of Theorem~\ref{thm: solvingset} asserts that every Cayley digraph $Y'$ on $G$ isomorphic to $X$ admits an isomorphism $f:X\to Y'$ for which $\langle G_L,f^{-1}G_Lf\rangle$ admits a complete sequence of block systems. Hence, for the present $Y$, there is such an $f:X\to Y$. Theorem~\ref{thm: main} then gives a color-preserving isomorphism $d=\prod_{p\mid |G|}d_p\colon X\to Y$ with $d_p\in\Sym(G_p)$. Since $d=\prod_{p\mid |G|}d_p\in D_G$ is a color-preserving Sylow coordinate isomorphism, Lemma~\ref{grp: coloursubstitution} replaces each factor $d_p$ by an element $s_p\in S_p$ and yields a color-preserving isomorphism $s=\prod_{p\mid |G|}s_p\colon X\to Y$.

Let $\prod_{p\mid |G|}S_p=\{(s_p)_{p\mid |G|}:s_p\in S_p\}$, where $(s_p)_{p\mid |G|}$ acts componentwise on $G=\prod_{p\mid |G|}G_p$, that is, $(s_p)_{p\mid |G|}(x_p)_{p\mid |G|}=(s_p(x_p))_{p\mid |G|}$. Then $s=(s_p)_{p\mid |G|}$ is a color-preserving isomorphism from $X$ to $Y$ belonging to this product. Since $Y$ was an arbitrary Cayley digraph on $G$ isomorphic to $X$, the same argument applies to every such $Y$. By the definition of a solving set, this product is a solving set for $X$.
\end{proof}

\section{Conclusion}\label{sec: conclusion}

This paper studies the isomorphism problem for Cayley digraphs on abelian groups of cube-free order. We prove that if two colored binary Cayley structures are connected by a color-preserving isomorphism $f$ for which $\langle G_L,f^{-1}G_Lf\rangle$ admits a complete sequence of block systems, then they are also connected by a color-preserving Sylow coordinate isomorphism. Consequently, for a fixed Cayley digraph $X$, the product of universal color-preserving solving sets, one for each Sylow subgroup, is a solving set for $X$, provided that every Cayley digraph isomorphic to $X$ is connected to $X$ by such an isomorphism. This yields a Sylow reduction and gives an affirmative answer to Dobson's open problem.

\backmatter

\section*{Statements and Declarations}

\subsection*{Funding}
Supported by the National Natural Science Foundation of China (Nos. 12261071, 12571019) and the Natural Science Foundation of Qinghai Province (No. \mbox{2025-ZJ-902T}).

\subsection*{Competing interests}
The authors declare that they have no conflicts of interest.

\end{document}